# Two-Stage Multi-Objective OPF for AC/DC Grids with VSC-HVDC: Incorporating Decisions Analysis into Optimization Process

Yang Li [a,b,*], Yahui Li [a], Guoqing Li [a], Dongbo Zhao [b], Chen Chen [b]

[a] School of Electrical Engineering, Northeast Electric Power University (NEEPU), Jilin 132012, China

[b] Center for Energy, Environmental, and Economic Systems Analysis (CEEESA), Argonne National Laboratory, IL 60439, United States

**ABSTRACT:** A two-stage solution approach for solving the problem of multi-objective optimal power flow (MOPF) is proposed for hybrid AC/DC grids with VSC-HVDC. First, a MOPF model for hybrid AC/DC grids is built to coordinate the economy, voltage deviation and environmental benefits. Then, a two-stage solution approach, incorporating decision analysis into optimization process, is presented to solve the model. The first stage of the approach is consisted of a multi-objective particle swarm optimization algorithm with a hybrid coding scheme employed to find multiple Pareto-optimal solutions. The second stage will have the obtained solutions clustered into different groups using fuzzy c-means (FCM) clustering, and then the 'best' compromise solutions are obtained by calculating the priority memberships of the solutions belonging to the same groups via grey relation projection (GRP) method. The novelty of this approach lies primarily in incorporating the FCM-GRP based decisions analysis technique into MOPF studies, thereby assisting decision makers to automatically identify the 'best' operation points. The effectiveness of the proposed approach is verified based on the test results of the IEEE 14- and 300- bus systems.

[*]Corresponding author. E-mail address: liyang@neepu.edu.cn (Yang Li).

## 1. Introduction

Optimal power flow (OPF) has always been regarded as one of the most important means for achieving optimal control and scheduling of the power systems [1, 2]. Traditionally, the main objective of OPF problem is to seek the minimum of economic costs, such as electricity generation costs or active power losses, by satisfying a set of various constraints [3]. With trends toward interconnections of AC/DC girds, a gradual deepening of electricity market reforms and large-scale integration of renewable energy sources, the system is operated approaching its stability limits. As a result, the conventional mono-objective OPF may not be adequate to analyze the increasingly complex and stressed power systems [4]. For this purpose, multi-objective OPF (MOPF) has become a hot topic in the field of OPF, since it adapts to the utilities' actual needs in coordinating multiple different weight- or even conflicting operational objectives [5-7].

As a new high voltage direct current (HVDC) transmission technology, voltage source converter based HVDC (VSC-HVDC) has been playing an increasingly important role in shaping the future of electric industry [8-10]. Compared with the conventional line-commutated converter based HVDC (LCC-HVDC) technology, VSC-HVDC has several significant advantages, such as enabling the integration of renewable energy sources, independent regulation of active and reactive powers, and capabilities of power supply for passive networks [11-13]. In addition, successful implementations of a series of major commercial VSC-HVDC projects have brought with it huge amounts of economic and social benefits [14]. Thereby, the OPF problem for AC/DC grids with VSC-HVDC is of paramount importance.

Numbers of recent attempts have been made to solve the mono-objective OPF problem for AC/DC grids with VSC-HVDC [15-20], but there has been little attention given to the MOPF issues of such systems. In [15], a second order cone programming has been used to solve the OPF problem for AC/DC systems with voltage source converter (VSC) technologies. In [16], the OPF problem considering the grid code constraints within VSC based multi-terminal DC (VSC-MTDC) is formulated as a master problem and its sub-problems by Benders' decomposition. In [17], a hybrid transmission grid architecture is presented to enable the efficient and globally optimal solution of the OPF problem. Reference [18] proposes an OPF model based on information-gap decision theory for handling the OPF problem with HVDC connected wind farms. Reference [19] presents a genetic algorithm (GA) to solve the OPF problem while maintaining an $N$-1 security constraint for grid integration of offshore wind farms with VSC-MTDC. In [20], an improved corrective security constrained

OPF formulation is proposed for meshed AC/DC grids with VSC-MTDC.

As is known, MOPF is a typical non-linear, non-convex, non-smooth, and high-dimension optimization problem. There is no optimal solution to make all the objectives to achieve optimal values simultaneously [4], and only a set of compromise solutions, called Pareto-optimal solutions, can be obtained. Unfortunately, it is difficult for decision makers to determine 'best' compromise solutions from all the Pareto optimals for two main reasons: (1) the size of the Pareto-optimal set is generally large and that decision vectors in the set contain different information representing preferences of decision makers; (2) different decision makers may have different preferences under the same operational condition, and even a same decision maker's preference may vary with time according to the actual operation requirements for his specific system. Therefore, how to automatically identify the 'best' compromise solutions becomes a challenging and urgent task for us to handle the MOPF issue.

*A. Contribution of This Paper*

This paper presents a two-stage solution approach to slove the problem of MOPF for AC/DC grids with VSC-HVDC. The use of multi-objective evolutionary algorithms (MOEAs) for MOPF has been deeply studied in [4-7]. However, these methods cannot be applied directly to such AC/DC grids, since VSC-HVDC has not been considered in their models. Moreover, all these methods cannot provide any decision support for decision makers to find the 'best' compromise solutions from the set of Pareto optimals. For this reason, the proposed approach introduces the decision analysis technique into the optimization process to overcome the disadvantages of existing MOPF methods lacking the ability to extract 'best' operating points.

The main contribution of this paper is two-fold. First, a MOPF model is built to coordinate the economy, voltage deviation, and environmental benefits for hybrid AC/DC grids. More importantly, the proposed two-stage method addresses how to employ the built MOPF model in detail. In the first stage, multi-objective particle swarm optimization (MOPSO) algorithm [21] with a hybrid coding scheme is used for finding well-spread Pareto-optimal solutions; in the second stage, the obtained solutions are clustered into different groups using fuzzy c-means (FCM) clustering [22], then the priority memberships of the solutions belonging to the same groups are calculated by the grey relation projection (GRP) method [23] to choose the 'best' compromise solutions.

*B. Organization of This Paper*

The remainder of this paper is organized as follows. First of all, VSC-HVDC steady-state model is briefly introduced. Then, a MOPF model for AC/DC Grids with VSC-HVDC is described in detail. Formulations of the proposed two-stage approach which incorporates FCM-GRP based decision analysis into optimization process are presented next. Application of the proposal is then demonstrated using two IEEE test systems, and finally the conclusions are made.

## 2. Steady-state Model of VSC-HVDC

### 2.1. Steady-state Power Characteristic

The simplified equivalent model of AC/DC system with VSC-HVDC is shown in Fig. 1.

In this model, $Z = R + jX$ represents the impendence between the AC and DC network, $\dot{U}_s = U_s \angle \delta_s$ is the AC bus voltage and $\dot{U}_c = U_c \angle \delta_c$ indicates the converter voltage, the current injected into converters from the AC grid is

$$\dot{I} = \frac{\dot{U}_s - \dot{U}_c}{R + jX} \tag{1}$$

The injected apparent power is given by

$$\dot{S}_s = P_s + jQ_s = \dot{U}_s \dot{I}^* \tag{2}$$

According to (2), the active power $P_s$ and the reactive power $Q_s$ injected from the AC network are obtained as

$$\begin{aligned} P_s &= U_s^2 G - U_s U_c [G\cos(\delta_s - \delta_c) + B\sin(\delta_s - \delta_c)] \\ Q_s &= -U_s^2 B - U_s U_c [G\sin(\delta_s - \delta_c) - B\cos(\delta_s - \delta_c)] \end{aligned} \tag{3}$$

where $G = 1/R$, $B = 1/X$.

The injected active power $P_c$ and the reactive power $Q_c$ from the AC grid into VSC are as follows.

$$\begin{aligned} P_c &= -U_c^2 G + U_s U_c [G\cos(\delta_s - \delta_c) + B\sin(\delta_s - \delta_c)] \\ Q_c &= U_c^2 B - U_s U_c [G\sin(\delta_s - \delta_c) + B\cos(\delta_s - \delta_c)] \end{aligned} \tag{4}$$

Due to the power injected into the DC grid $P_{dc}$ is equal to the active power $P_c$, the following equation can be derived.

$$P_{dc} = -U_c^2 G + U_s U_c [G\cos(\delta_s - \delta_c) + B\sin(\delta_s - \delta_c)] \tag{5}$$

**2.2. Control strategies of VSC-HVDC**

Basic control strategies of VSC-HVDC can be divided into four classes [24]: 1) constant DC voltage and reactive power control; 2) constant DC voltage and AC voltage control; 3) constant active and reactive power control (constant *PQ*-control); 4) constant active power and AC voltage control (constant *PV*-control).

As far as VSC-MTDC, only one terminal must adopt constant the DC voltage control to maintain the power balance, while other terminals may employ constant *PQ*- or *PV*-control. In addition, droop control has been becoming an appealing control strategy for VSC-MTDC.

## 3. Problem Formulation

In this section, a MOPF model is introduced for coordinating economy, voltage deviation, and environmental benefits of AC/DC system with VSC-HVDC.

### 3.1. Objective Functions

#### 3.1.1. System Active Power Losses

The system active power losses of AC/DC grids is given by

$$O = \sum_{i=1}^{N} U_i \sum_{j=1}^{N} U_j (G_{ij} \cos \theta_{ij} + B_{ij} \sin \theta_{ij}) + P_{DC\_loss} \tag{6}$$

where $O$ is the total active power losses of AC/DC grids, $N$ is the number of buses in the AC grid, $U_i$ and $U_j$ are respectively the voltage amplitudes of bus $i$ and bus $j$; $G_{ij}$, $B_{ij}$ and $\theta_{ij}$ are respectively the conductance, susceptance and phase-angle difference between buses $i$ and $j$. $P_{DC\_loss}$ is the power losses in the DC grid, comprising power losses of converter stations and DC lines. As the main part of $P_{DC\_loss}$, the converter station loss $P_{DC\_loss}$ is determined by the converter current $I_c$ [24]

$$\begin{cases} P_{con\_loss} = a + b \cdot I_c + c \cdot I_c^2 \\ I_c = \dfrac{\sqrt{P_c^2 + Q_c^2}}{\sqrt{3}U_c} \end{cases} \quad (7)$$

where $a$, $b$ and $c$ are the per unit coefficients. The losses of DC lines can be computed according to voltages and impedances.

### 3.1.2. Emissions of Polluting Gases

The utilized model of polluting gas emissions is [25]

$$E = \sum_{i=1}^{N_G}(\alpha_i P_{G,i}^2 + \beta_i P_{G,i} + \gamma_i) \quad (8)$$

where $E$ is the amount of polluting gas emissions, $N_G$ is the number of generators, $P_{G,i}$ is the active power output of the $i$th generator, $\alpha_i$, $\beta_i$ and $\gamma_i$ are respectively the factors of polluting gases emissions of generator $i$.

### 3.1.3. Voltage Deviation Index

For improving the voltage quality and operational security of the AC/DC grid, the voltage deviation index $V_{de}$ is used as the objective function.

$$V_{de} = \sum_{j=1}^{N}(U_j - U_{set,j})^2 + \sum_{k=1}^{N_{dc}}(U_{dc,k} - U_{set,dc,k})^2 \quad (9)$$

where $U_j$ is the voltage of AC bus $j$, and $U_{set,j}$ indicates the predefined value of $U_j$; $U_{dc,k}$ is the voltage of DC bus $k$, and $U_{set,dc,k}$ indicates the predefined value of $U_{dc,k}$.

## 3.2. Constraints

### 3.2.1. Equality Constraints

The equality constraints of the AC grid are

$$\begin{cases} P_{gi} - P_{di} - U_i \sum_{j \in i} U_j \left( G_{ij} \sin \theta_{ij} + B_{ij} \cos \theta_{ij} \right) = 0 \\ Q_{gi} - Q_{di} - U_i \sum_{j \in i} U_j \left( G_{ij} \sin \theta_{ij} - B_{ij} \cos \theta_{ij} \right) = 0 \end{cases} \quad (10)$$

where $P_{gi}$ and $Q_{gi}$ are the active and reactive power inputs of bus $i$; $P_{di}$ and $Q_{di}$ are the active and reactive loads of bus $i$.

Besides equation (3)-(5), the equality constraints of the DC grid includes

$$\sum_{i=1}^{N_{dc}} P_{s,i} + P_{DC\_loss} = 0 \quad (11)$$

where $P_{s,i}$ is the active power of bus $i$ in the AC grid, $N_{dc}$ is the number of buses in the DC grid.

### 3.2.2. Inequality Constraints

The inequality constraints in the AC grid can be written as

$$\begin{cases} P_{G,i}^{\min} \leq P_{G,i} \leq P_{G,i}^{\max}, & i = 1, ..., N_G \\ Q_{G,i}^{\min} \leq Q_{G,i} \leq Q_{G,i}^{\max}, & i = 1, ..., N_G \\ U_i^{\min} \leq U_i \leq U_i^{\max}, & i = 1, ..., N \\ T_i^{\min} \leq T_i \leq T_i^{\max}, & i = 1, ..., N_T \\ Q_{C,i}^{\min} \leq Q_{C,i} \leq Q_{C,i}^{\max}, & i = 1, ..., N_C \end{cases} \quad (12)$$

where $P_G$ and $Q_G$ are the generator active and reactive power output; $U$ is the node voltage; $T$ is the transformer tap ratio; $Q_C$ is the capacity of a reactive power compensation equipment; the superscript "max" and subscript "min" are respectively the corresponding upper and lower limits of the physical quantities; $N_T$ and $N_C$ are the number of transformers and reactive power compensation equipment.

The related inequality constraints of the converter stations in the DC grid [24] are

$$r_{\min}^2 \leq (P_s - P_0)^2 + (Q_s - Q_0)^2 \leq r_{\max}^2 \quad (13)$$

The equation is constraint range of active and reactive powers, $S_0(P_0, Q_0)$ is the center of the circle formed by the PQ-capability of each converter, $r_{max}$ and $r_{min}$ are respectively the upper and lower limits of the radius $r$ of the circle [24].

The inequality constraints in the DC grid can be written as

$$\begin{cases} P_{s,i}^{min} \leq P_{s,i} \leq P_{s,i}^{max}, & i=1,...,N_{dc} \\ Q_{s,i}^{min} \leq Q_{s,i} \leq Q_{s,i}^{max}, & i=1,...,N_{dc} \\ U_{dc,i}^{min} \leq U_{dc,i} \leq U_{dc,i}^{max}, & i=1,...,N_{dc} \\ I_{dc,ij}^{min} \leq I_{dc,ij} \leq I_{dc,ij}^{max}, & i,j=1,...,N_{dc} \end{cases} \quad (14)$$

where $P_s^{max}$ and $P_s^{min}$ are respectively the upper and lower limits of the active power $P_s$ of each terminal in the DC grid; $Q_s^{max}$ and $Q_s^{min}$ are respectively the upper and lower limits of the reactive power $Q_s$ of each terminal in the DC grid; $U_{dc}^{max}$ and $U_{dc}^{min}$ are respectively the upper and lower limits of the DC voltage $U_{dc}$ of each terminal in the DC grid; $I_{dc,ij}$ is the current of DC branches from bus $i$ to bus $j$, $I_{dc,ij}^{max}$ and $I_{dc,ij}^{min}$ are the upper and lower limits, respectively.

## 4. Proposed Two-stage Solution Approach

### 4.1. Overall Solution Framework

The overall framework for solving the built MOPF model incorporates two stages: multi-objective optimization and decision analysis, as shown in Fig. 2.

The MOPSO proposed by Coello Coello [21] is a powerful optimization algorithm for handling problems including both continuous and discrete variables. It has significant advantages over other MOEAs, such as fast converge speed and strong optimization ability, and has been widely used in various multi-objective problems in engineering [26, 27]. Here, it is used to find the Pareto-optimal solutions, as shown in Fig. 3.

The main steps of the optimization process are listed as follows.

Step 1: Input the initial variables, comprising: 1) AC/DC system parameters (the data of buses, branches, loads and generators, the number and control strategies of VSCs); 2) MOPSO parameters (population size, repository size, maximum iteration number, mutation rate and so on); 3) the ranges of related variables and the steps of discrete variables.

Step 2: Initialize population. According to the characteristics of controlled variables, a hybrid coding scheme is used to facilitate the optimization since continuous and discrete variables require different coding schemes. The continuous variables comprise the active power output $P_G$ of a generator (except the balancing machine), the generator terminal voltage $U_G$, $P_s$, $Q_s$, and $U_{dc}$, while the discrete variables comprise $T$ and $Q_C$. And thereby, a particle corresponding to an OPF solution is coded as

$$\left( \underbrace{P_{G,2},\cdots,P_{G,N_G},U_{G,1},\cdots,U_{G,N_G},P_{s,1},\cdots,P_{s,N_{dc}},Q_{s,1},\cdots,Q_{s,N_{dc}},U_{dc,1},\cdots,U_{dc,N_{dc}}}_{\text{continuous variables}}, \underbrace{T_1,\cdots,T_{N_T},Q_{C,1}\cdots,Q_{C,N_C}}_{\text{discrete variables}} \right).$$

In this way, the initial population is generated uniformly at random throughout the entire feasible search space. And then, set the iteration counter $I_{ite}$ to 0.

Step 3: Via the alternating iterative method proposed by Beerten [28], the AC/DC system power flow is calculated with the constraints considered, and the values of objective functions $f=(O,E,V_{de})$ for each particle are obtained.

Step 4: Add non-dominated solutions into the external archive by comparing with each particle. Determine the initial optimal position of each particle, called *pbest*, and global optimal position of all particles, called *gbest*.

Step 5: Calculate and adjust the positions and velocities of all particles, ensuring the particles flying within the search space according to the optimization mechanism of the algorithm. The current velocities are determined on *pbest*, *gbest*, the previous positions and velocities of particles, the current positions is the sum of previous positions and current velocities [21].

Step 6: Compute the values of objective functions $f$ after updating AC/DC gird data with current positions, adjust *pbest* for each particle in the new particle swarm (if the current position is dominated by *pbest*, keep *pbest* unchanged; otherwise, replace *pbest*).

Step 7: Archive maintenance based on the non-dominated solutions, and select *gbest* for all particles in the population.

Step 8: Judge the termination criteria. If $I_{ite}$ exceeds the preassigned maximum iteration number, execute Step 9; otherwise, increase $I_{ite}$ by 1 and go to Step 5 to repeat the above process.

Step 9: Output the obtained Pareto-optimal solutions.

### 4.3. Second Stage: FCM-GRP Based Decision Analysis

In this section, a novel FCM-GRP based decision analysis technique is described to identify the 'best' compromise solutions from all the Pareto-optimal solutions.

#### 4.3.1. FCM Clustering

FCM proposed by Bezdek is a classical unsupervised clustering algorithm based on solving the following problem [22].

$$\min \ F_t(W,U,V) = \sum_{i=1}^{N_p}\sum_{j=1}^{c} \mu_{ij}^t \|w_i - v_j\|^2$$

$$s.t. \quad \sum_{j=1}^{c} \mu_{ij} = 1$$

(15)

where $F_t$ is the loss function defined by membership functions, $W = \{w_1, w_2, \cdots, w_i, \cdots, w_{N_p}\}$ is the input vector of $N_p$ Pareto-optimal solutions to be clustered, $U = \{u_1, u_2, \cdots, u_c\}$ and $V = \{v_1, v_2, \cdots, v_j, \cdots, v_c\}$ are respectively the output vectors of membership degrees and cluster centers, $\mu_{ij}\ (\mu_{ij} \in [0,1])$ is the membership degree of $x_i$ to the cluster $v_j$, $c$ is the pre-defined cluster number, and $t\ (t \in [1, \infty])$ is a parameter controlling the fuzziness of the clustering process.

The principle of FCM is to minimize the loss function by iteratively updating membership functions and all clustering centers. The process iteratively repeats until a pre-given tolerance is satisfied or the generation exceeds the maximum number of iterations. To represent the preferences of decision makers over three objective functions, the clustering number is set to 3 in this work. In this way, FCM provides the cluster centroids as representative data points of the Pareto-optimal solutions for the problem and maps each data point to a cluster.

**4.3.2. Grey Relation Projection**

GRP theory is a powerful tool for analyzing the relationship between sequences with grey information and has been successfully applied in a variety of fields [23, 29]. Since all the three optimized objective functions belong to "benefit-type" indicators according to their characteristic, the projection $V_l$ of each scheme on the ideal reference scheme is given by

$$V_l^{+(-)} = \sum_{k=1}^{t} \gamma_{lk}^{+(-)} \frac{w_k^2}{\sqrt{\sum_{k=1}^{t}(w_k)^2}} \qquad (16)$$

where superscript "+" denotes positive scheme, superscript "−" is negative scheme, $t$ is the number of indication, $\gamma_{lk}$ is the grey relation coefficient between the $k$th indication in $l$th scheme, $w_k$ is the weight of each indication in the scheme. For the convenience of description the weights of three objective functions are equal in this study, but they can be arbitrarily adjusted by decision makers on the basis of their preferences. The priority membership $d$ is

$$d_l = \frac{(V_0 - V_l^-)^2}{(V_0 - V_l^-)^2 + (V_0 - V_l^+)^2}, \quad 0 \leq d_l \leq 1 \qquad (17)$$

where $V_0$ is equal to $V_l$ when $\gamma$ takes 1. The equation shows that the higher the priority membership is, the better the scheme will be. Therefore, the results with the highest priority membership are chosen as the 'best' compromise solutions.

**5. Case Study**

To verify the effectiveness of the proposed approach, a modified IEEE 14-bus system is firstly used as the test system. This system includes 5 generators, 18 branches, 11 loads and modular multilevel converters [15].

All of the programs are implemented in the MATLAB environment on a desktop computer with the master frequency 3.40 GHz and main memory 4 GB. Note that, for ease of description, only the computation time in the first stage, rather than the total simulation time, is provided and discussed hereinafter. The reason for this is that the computation time in the second stage can be negligible (only about

0.1 s), which is far less than that in the first stage.

## 5.1. Parameter Setting

### 5.1.1. Parameters of Hybrid AC/DC Grids

The base capacity of the system is 100 MVA. The voltage of each bus ranges from 0.94 to 1.06 pu; the tap ratio of each transformer is in the range from 0.9 to 1.1, and the step size is 0.0125; the capacity of reactive power compensation equipment installed at bus 9 ranges from 0 to 0.25 pu, and the step size is 0.01 pu; the $P_s$ and $Q_s$ are in the range of [-1.0 pu, 1.0 pu], and the DC voltage ranges from 0.9 to 1.1 pu.

### 5.1.2. Parameters of MOPSO

The main parameters of MOPSO are listed in Table 1.

## 5.2. IEEE 14-bus System with Two-terminal DC Network

In this case, the IEEE 14-bus system is embedded into a two-terminal DC network, and the branch 4-5 is modified to a VSC-HVDC link. The modified system is shown in Fig. 4.

The initial DC bus parameters are listed in Table 2, where $B_{VSC}$ represents the bus number connected to a VSC.

Taking the modified system as the test case, the optimization effects of proposed approach are analyzed. The adopted control strategy of VSC$_1$ is constant DC voltage and reactive power control, and the control strategy of VSC$_2$ is constant *PQ*-control. It should be noted that the proposal is also applicable to other control strategies besides the above-mentioned ones.

### 5.2.1. OPF with two Optimized Objectives

In order to examine the effectiveness of the proposal, the well-known non-dominated sorting genetic algorithm (NSGA-II) and MOPSO are respectively employed to seek the Pareto-optimal solutions of the MOPF problem, taking the system active power losses $O$ and emissions of polluting gases $E$ as the objective functions. For NSGA-II, the population size and the maximum iteration number are set as same as those of MOPSO, the crossover probability $pc$ is 0.9 and mutation probability $pm$ is $1/L_c$ ($L_c$ is the length of a chromosome). Fig. 5 gives the changing proportion of the Pareto non-dominated solutions to the total population with the increase of iterations.

As shown in Fig. 5, the convergence rate of the MOPSO algorithm is superior to that of NSGA-II. When stably reaching the Pareto frontier, NSGA-II needs to repeatedly iterate 25 epochs, while MOPSO only requires 15 iterations.

When the predefined maximum iterations are reached, the distribution of Pareto-optimal solutions of NSGA-II and MOPSO in the objective function space are shown in Fig. 6, and the extreme solutions are listed in Table 3.

From Fig. 6 and Table 3, it can be seen that:

The two objectives, minimize active power loss $O$ and emissions of polluting gases $E$, are conflicting to each other. In other words, seeking the minimum of $O$ will inevitably lead to the increase of $E$ at the expense; and vice versa.

Considering the inherent randomness of MOEAs, all experiments using MOPSO and NSGA-II have been repeatedly run 30 times. Table 4 gives the required maximum, minimum and average number of iterations when the algorithms stably reach their own Pareto frontiers.

Table 4 shows that MOPSO has a faster convergence rate than NSGA-II when reaching Pareto frontiers.

Table 5 gives the maximum, minimum and average computation times of the two algorithms in the 30 runs.

From Table 5, it can be clearly seen that the computation times of MOPSO are less than that of NSGA-II in all the statistical indicators. The results demonstrate that MOPSO excels over NSGA-II in terms of computation efficiency.

In recent years, the quantitative assessment of the performance of an MOEA has been attracting an increasing attention. However, up to now, there is still no general agreement on the criteria to evaluate the Pareto frontiers obtained by different MOEAs. In general, an effective performance measure of an MOEA musts satisfy the following basic criteria [30-32]. a) convergence: the computed Pareto front should be as close as possible to the real Pareto' front; b) distribution of solutions: the Pareto solutions should be evenly distributed with a good diversity along the Pareto front.

In order to properly evaluate the solution qualities of NSGA-II and MOPSO, the following two quantitative metrics are thereby taken into account corresponding to the above criteria in this work.

**1) Generational distance**

The generational distance (GD) is a well-known indicator to estimate how far elements are in the set of non-dominated vectors found so far from those in the Pareto-optimal set [21, 30-31]. It is defined as follows,

$$GD = \frac{\sqrt{\sum_{i=1}^{N_s} D_i^2}}{N_s} \quad (18)$$

where $N_s$ is the number of non-dominated solutions, $D_i$ denotes Euclidean distance between each of these and the nearest member of the Pareto optimal set. The GD can represent the convergence metric, and the less value of GD shows the better condition of the related Pareto-optimal solutions set.

**2) Spacing**

The spacing (SP) is another powerful metric to measure the distribution of the Pareto-optimal solutions, and it is given in the following formulation [21, 30-32]:

$$SP = \sqrt{\frac{1}{N_s - 1} \sum_{i=1}^{N_s} (\bar{D} - D_i)^2} \quad (19)$$

where $\bar{D}$ is the mean of $D_i$. Similar to the GD metric, the less value of SP means the better distribution of the Pareto-optimal solutions.

Table 6 shows the statistical results of the metrics for two algorithms in 30 runs. From the above table, it can be observed that all indicator values of the MOPSO are superior or equal to that of NSGA-II. The results demonstrate that the MOPSO outperforms the NSGA-II in solution quality, i.e. the Pareto frontiers obtained by the MOPSO algorithm achieves better solution quality than NSGA-II. Therefore, the conclusions can be drawn: a) MOPSO is an effective tool to seek the Pareto-optimal solutions for the built MOPF model; b) In terms of both computation efficiency and solution quality, the MOPSO is superior to the commonly used NSGA-II.

**5.2.2. OPF with Three Optimized Objectives**

MOPSO algorithm is first utilized to obtain the Pareto-optimal solutions, and the average computation time is 26.21 s by using MOPSO in 30 runs. Without loss of generality, one of the results is taken as an example for further analysis and its distribution in the objective function space is shown in Fig. 7.

As shown in Fig. 7, it is clear that the proposal can produce nearly complete and uniform Pareto optimals. And thereby, the conclusion can be drawn that the multi-objectives in MOPF can be coordinated by using the proposed approach.

The extreme solutions corresponding to the different objective functions are shown in Table 7.

From Table 7, it can be observed that the each extreme solutions respectively reach the minimum values of $O$, $E$ and $V_{de}$. Note that, as far as mono-OPF, the obtained three extreme solutions are optimal; but for MOPF, they are non-inferior solutions.

After obtaining the Pareto-optimal solutions, FCM clustering is employed to cluster the optimals into three groups marked with different colors, as shown in Fig. 8. Note that in this work the red, green and blue color respectively represent decision makers' preference to the objective functions $O$, $E$ and $V_{de}$.

Then, GRP is utilized to calculate the priority memberships of the Pareto optimals belonging to the same groups, and the solutions with the highest memberships are chosen as the 'best' compromise solutions, as shown in Table 8.

When a decision-maker tends to put more emphasis on the economy index $O$, the compromise solution I might be a preferable choice; when the environmental benefits $E$ is taken as the first objective, then the compromise solution II is a more valuable option than the others; when placing top priority on the role of voltage deviation index $V_{de}$, the compromise solution III might be considered the best one. In this way, the proposed approach can give the 'best' compromise solution from all the Pareto-optimal solutions by incorporating FCM-GRP based decision analysis into the optimization process, helping decision makers to automatically determine the 'best' operating.

Taking compromise solution I in Table 8 as an example, some variables before and after optimization are shown in Tables 9-11.

Tables 9-11 demonstrate that the power flow distribution of the system is more reasonable after optimization, embodying that all the objective functions are improved to different extents.

### 5.3. IEEE 14-bus System with Three-terminal DC Network

In order to further evaluate the effectiveness of the proposal for AC/DC grids with VSC-MTDC, it is tested using a modified IEEE 14-bus system with three-terminal DC network, as shown in Fig. 9. The $VSC_1$, $VSC_2$, and $VSC_3$ are respectively connected with buses 2, 4 and 5 in the AC part, and the initial DC bus parameters are shown in Table 12.

#### 5.3.1. 'Best' Compromise Solutions

Here, the used control strategy of $VSC_1$ is constant DC voltage and reactive power control, and the control strategies of both $VSC_2$ and $VSC_3$ are constant *PQ*-control. Note that, other control strategies, e.g. droop control, are also applicable to the proposed approach.

The average computation time in 30 runs is 26.77 s by using MOPSO in the first stage. Similar in Section 5.2.2, one of the results is utilized in the second stage. The distribution of this result in the objective function space after FCM clustering is shown in Fig. 10, and the recommended 'best' compromise solutions are listed in Table 13.

From Fig. 10, it shows that well-distributed Pareto-optimal solutions can be obtained in this case. Table 13 demonstrates that the 'best' compromise solutions can be automatically identified as well. Consequently, the proposed method is readily applicable to the system with VSC-MTDC.

### 5.3.2. Comparative Results Before and After Optimization

To demonstrate the effectiveness of the proposal, comparative results before and after optimization are discussed by taking the recommended solution I in Table 11 as an example. Some variables before and after optimization are shown in Tables 14-16.

According to the Tables 14-16 the power flow distribution after optimization is more reasonable.

### 5.4. Influence of Different Control Strategies and Number of VSCs

To verify the influence of different control strategies and the number of VSCs, three typical operation conditions are considered.

Case1——both the system configuration and control strategies are the same as the AC/DC system in section 5.2;

Case2——the used system configuration is the same as Case1, but the control strategies of two VSCs are exchanged;

Case3——both the system configuration and control strategies are the same as the AC/DC system in section 5.3.

After the optimization process in the first stage, the solutions with the highest priority memberships in the three conditions are employed as the reference solutions, and the corresponding operation variables are shown in Table 17. Here, the 'maximum' and 'minimum' in the table respectively denote the upper- and lower- limits of the variables.

Table 17 shows that the used control strategies have a certain impact on optimization of power flow. Specially speaking, $O$ and $V_{de}$ in Case1 are superior to those in Case2, but $E$ in Case2 is better than that in Case1. In addition, Table 17 shows that all objective functions in Case3 are better than those in Case1 and Case2. The reason for this is that the regulating ability of DC grid is correspondingly enhanced with the increase of the number of VSCs, so that the distribution of power flow tends to be more reasonable.

## 6. Application to the IEEE 300-bus System

For the purpose of examining the applicability of the proposal to a larger-scale system, it is further tested on the power system of IEEE 300-bus system. This system is made up of 69 generators, 407 branches, 201 loads, and some series compensation devices and static var compensators. Meanwhile, four VSCs are respectively connected with buses 99, 105, 107 and 110 in the 300-bus system. The original branches between the four buses are modified to VSC-HVDC links. The control strategy of VSC connected to bus 105 is constant DC voltage and reactive power control, and other VSCs adopt constant *PQ*-control.

The average computation time is 115.34 s by using MOPSO in 30 runs. Taking the obtained Pareto-optimal solution set in one run as example, the distribution of the solutions after clustering is shown in Fig. 11.

Fig. 11 shows that the Pareto-optimals with a good distribution can also be found. Table 18 gives the 'best' compromise solutions via GRP method.

Taking compromise solution I as an example, the comparative results before and after optimization are listed in Table 19.

As is shown in Table 19, the proposal has achieved good optimization results. Specifically speaking, the network loss rate is reduced from 1.31% to 1.17%; the amount of polluting gas emissions is reduced from 340453 lb/h to 316016 lb/h, and the voltage deviation index is decreased from 0.2465 to 0.1889.

The above results suggest that the proposal is also effective for IEEE 300-bus system, and thereby, the applicability of the proposed MOPF model and solution approach to the larger-scale power system is verified.

## 7. Conclusion

During the past few years, the emerging VSC-HVDC technology has become an increasingly appealing option for bulk power transmission. To solve the problem of MOPF for AC/DC grids with VSC-HVDC, this paper has presented a two-stage solution approach by incorporating the FCM-GRP based decision analysis techniques into the optimization process. Studies carried out on IEEE 14- and 300- bus systems prove that the proposed approach is effective in handling this issue. The built MOPF model can coordinate the economy, voltage deviation and environmental benefits in a unified manner, thereby adapting to the utilities' actual needs of coordinating multiple operational objectives. More importantly, the proposed method not only provides the well-distributed Pareto-optimal solutions via multi-objective optimization, but also gives the 'best' compromise solutions automatically with the use of decision analysis. The method will find extensive potential applications in the intelligent planning and operation of smart grids.

Future work will focus on improving the efficiency and practicality of the proposed approach for physical applications. The efficiency of the proposed method can be further elevated through the use of parallel computing techniques when applied in large-scale real systems. In addition, this work will be extended to a larger framework of security-constrained MOPF.

## Acknowledgements

This work was supported by the China Scholarship Council (CSC), the U.S. Department of Energy (DOE)'s Office of Electricity Delivery and Energy Reliability – Advanced Grid Modeling (AGM) Program, the National Key Research, Development Program of China under Grant No. 2017YFB0902401, and the National Natural Science Foundation of China under Grant No. 51677023.

**References**

[1] M. Ghasemi, S. Ghavidel, E. Akbari, and A. A. Vahed. Solving non-linear, non-smooth and non-convex optimal power flow problems using chaotic invasive weed optimization algorithms based on chaos. Energy 2014; 73: 340-353.

[2] C. Coffrin, H. L. Hijazi, P. Van Hentenryck. The QC relaxation: A theoretical and computational study on optimal power flow. IEEE Trans Power Syst 2016; 31(4): 3008-3018.

[3] R. A. Abarghooee, T. Niknam, M. Malekpour, F. Bavafa, M. Kaji. Optimal power flow based TU/CHP/PV/WPP coordination in view of wind speed, solar irradiance and load correlations. Energ Convers Manage 2015; 96: 131-145.

[4] M. Ghasemi, S. Ghavidel, M. M. Ghanbarian, M. Gharibzadeh, and A. A. Vahed. Multi-objective optimal power flow considering the cost, emission, voltage deviation and power losses using multi-objective modified imperialist competitive algorithm. Energy 2014; 78: 276-289.

[5] M. Pourakbari-Kasmaei, M. J. Rider, and J. R. Mantovani. An unequivocal normalization-based paradigm to solve dynamic economic and emission active-reactive OPF (optimal power flow). Energy 2014; 73: 554-566.

[6] X. Yuan, B. Zhang, P. Wang, J. Liang, Y. Yuan, Y. Huang, and X. Lei. Multi-objective optimal power flow based on improved strength Pareto

evolutionary algorithm. Energy 2017; 122: 70-82.

[7] S. Shargh, B. Mohammadi-ivatloo, H. Seyedi, and M. Abapour. Probabilistic multi-objective optimal power flow considering correlated wind power and load uncertainties. Renew Energ 2016; 94: 10-21.

[8] L. Q. Liu, and C. X. Liu. VSCs-HVDC may improve the electrical grid architecture in future world. Renew Sust Energy Rev 2016; 62: 1162-1170.

[9] A. Egea-Àlvarez, M. Aragüés-Peñalba, E. Prieto-Araujo, and O. Gomis-Bellmunt. Power reduction coordinated scheme for wind power plants connected with VSC-HVDC. Renew Energ 2017; 107: 1-13.

[10] O. A. Urquidez, and L. Xie. Smart targeted planning of VSC-based embedded HVDC via line shadow price weighting. IEEE Trans Smart Grid 2015; 6(1): 431-440.

[11] F. D. Bianchi, J. L. Domínguez-García, and O. Gomis-Bellmunt. Control of multi-terminal HVDC networks towards wind power integration: A review. Renew Sust Energy Rev 2016; 55: 1055-1068.

[12] A. Korompili, Q. Wu, and H. Zhao. Review of VSC HVDC connection for offshore wind power integration. Renew Sust Energy Rev 2016; 59, 1405-1414.

[13] M. A. Abdelwahed, and E. F. El-Saadany. Power sharing control strategy of multiterminal VSC-HVDC transmission systems utilizing adaptive voltage droop. IEEE Trans Sustain Energ 2017; 8(2): 605-615.

[14] J. Wu, Z. X. Wang, L. Xu, and G. Q. Wang. Key technologies of VSC-HVDC

and its application on offshore wind farm in China. Renew Sust Energy Rev 2014; 36, 247-255.

[15] M. Baradar, M. R. Hesamzadeh, and M. Ghandhari. Second-order cone programming for optimal power flow in VSC-type AC-DC grids. IEEE Trans Power Syst 2013; 28(4): 4282-4291.

[16] H. Y. Kim, and M. K. Kim. Optimal generation rescheduling for meshed AC/HIS grids with multi-terminal voltage source converter high voltage direct current and battery energy storage system. Energy 2017; 119: 309-321.

[17] M. Hotz, and W. Utschick. A hybrid transmission grid architecture enabling efficient optimal power flow. IEEE Trans Power Syst 2016; 31(6): 4504-4516.

[18] A. Rabiee, A. Soroudi, and A. Keane. Information gap decision theory based OPF with HVDC connected wind farms. IEEE Trans Power Syst 2015; 30(6): 3396-3406.

[19] R. T. Pinto, S. F. Rodrigues, E. Wiggelinkhuizen, R. Scherrer, P. Bauer, and J. Pierik. Operation and power flow control of multi-terminal DC networks for grid integration of offshore wind farms using genetic algorithms. Energies 2012; 6(1): 1-26.

[20] J. Cao, W. Du, and H. F. Wang. An improved corrective security constrained OPF for meshed AC/DC grids with multi-terminal VSC-HVDC. IEEE Trans Power Syst 2016; 31(1): 485-495.

[21] C. A. C. Coello, G. T. Pulido, and M. S. Lechuga. Handling multiple

objectives with particle swarm optimization. IEEE Trans Evolut Comput 2004; 8(3): 256-279.

[22] O. P. Mahela, and A. G. Shaik. Power quality recognition in distribution system with solar energy penetration using S-transform and Fuzzy C-means clustering. Renew Energ 2017; 106:37-51.

[23] X. Zhang, F. Jin, and P. Liu. A grey relational projection method for multi-attribute decision making based on intuitionistic trapezoidal fuzzy number. Appl Math Model 2013; 37(5): 3467-3477.

[24] J. Beerten, S. Cole, and R. Belmans. Generalized steady-state VSC MTDC model for sequential AC/DC power flow algorithms. IEEE Trans Power Syst 2012; 27(2): 821-829.

[25] A. M. Shaheen, R. A. El-Sehiemy, and S. M. Farrag. Solving multi-objective optimal power flow problem via forced initialised differential evolution algorithm. IET Gener Transm Dis 2016; 10(7): 1634-1647.

[26] J. Clarke, and J. T. McLeskey. Multi-objective particle swarm optimization of binary geothermal power plants. Appl Energ 2015; 138: 302-314.

[27] S. Meo, A. Zohoori, and A. Vahedi. Optimal design of permanent magnet flux switching generator for wind applications via artificial neural network and multi-objective particle swarm optimization hybrid approach. Energ Convers Manage 2016; 110: 230-239.

[28] J. Beerten, and R. Belmans. Development of an open source power flow

software for HVDC grids and hybrid AC/DC systems: MatACDC. IET Gener Transm Dis 2015; 9(10): 966-974.

[29] H. C. Liu, J. X. You, X. J. Fan, and Q. L. Lin. Failure mode and effects analysis using D numbers and grey relational projection method. Expert Syst Appl 2014; 41(10): 4670-4679.

[30] M. Q. Li, S. X. Yang, and X. H. Liu. Diversity comparison of Pareto front approximations in many-objective optimization. IEEE Trans Cybernetics, 2014; 44(12): 2568-2584.

[31] S. W. Jiang, Y. S. Ong, J. Zhang, and L. Feng. Consistencies and contradictions of performance metrics in multiobjective optimization. IEEE Trans Cybernetics, 2014; 44(12): 2391-2404.

[32] H. Lu, M. M. Zhang, Z. M. Fei, and K. F. Mao. Multi-objective energy consumption scheduling in smart grid based on Tchebycheff decomposition. IEEE Trans Smart Grid, 2015; 6(6): 2869-2883.

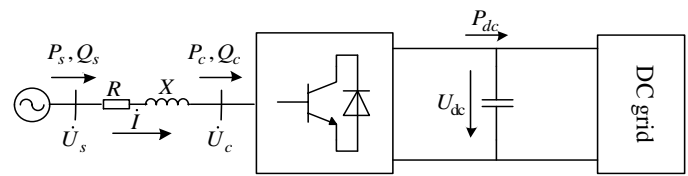

**Fig. 1.    Simplified equivalent model of AC/DC system**

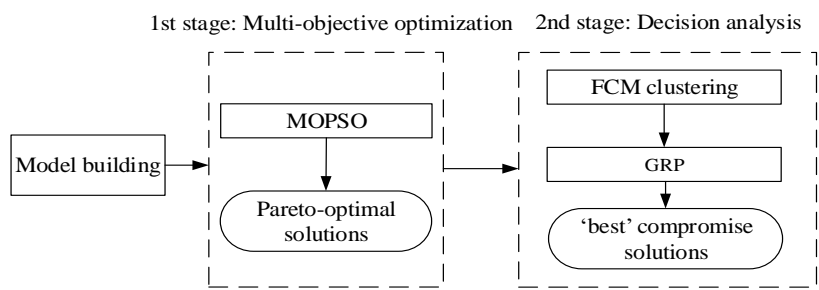

Fig. 2. Solution framework

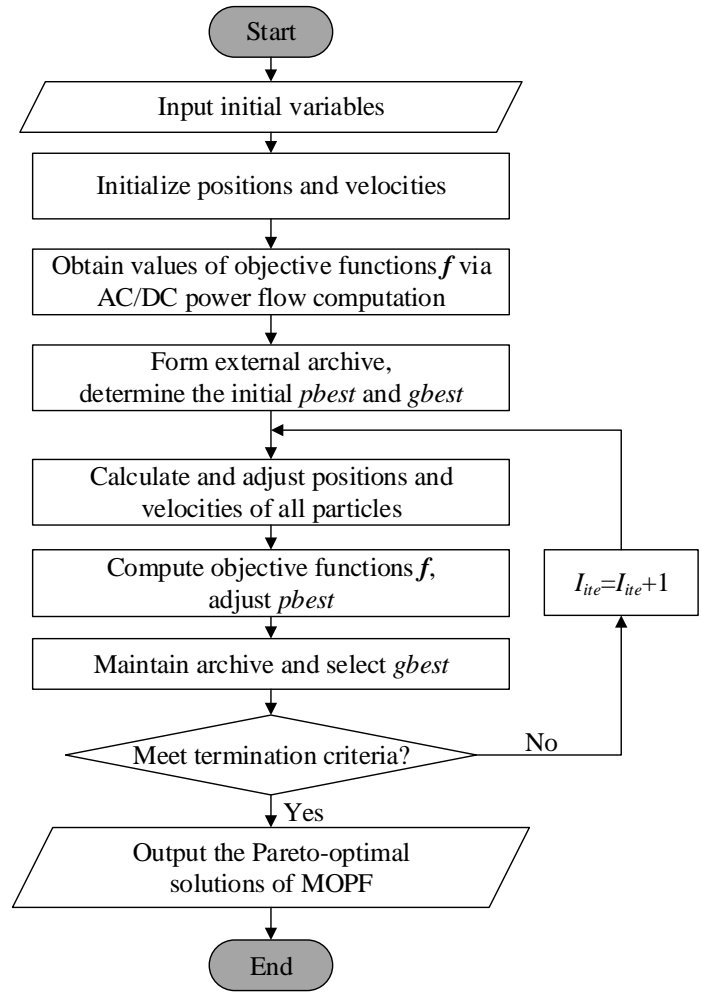

Fig. 3. Flowchart of optimization scheme using MOPSO

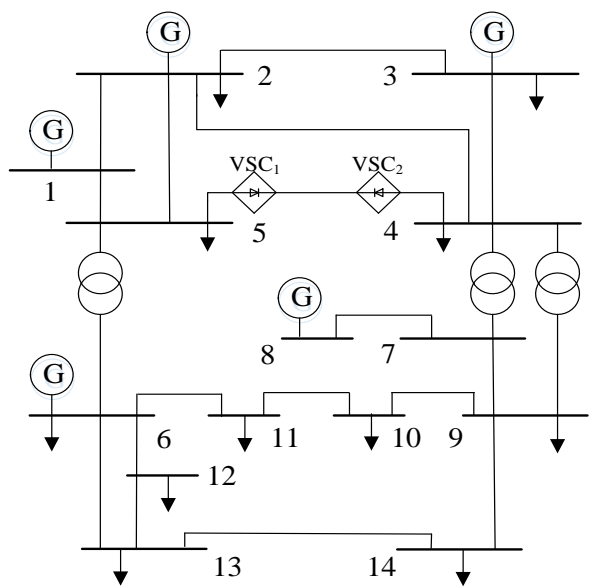

Fig. 4. IEEE 14-bus system with two-terminal DC network

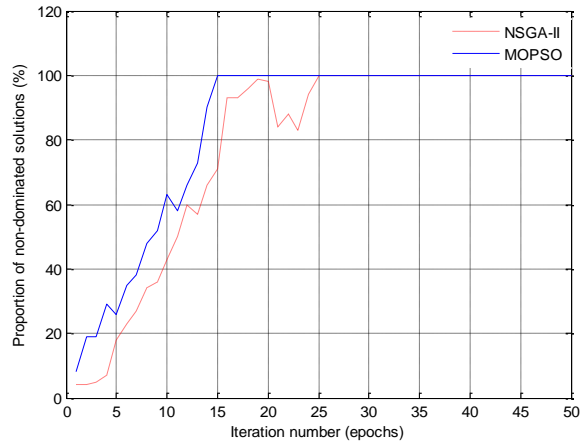

**Fig. 5. Proportion of Pareto non-dominated solutions**

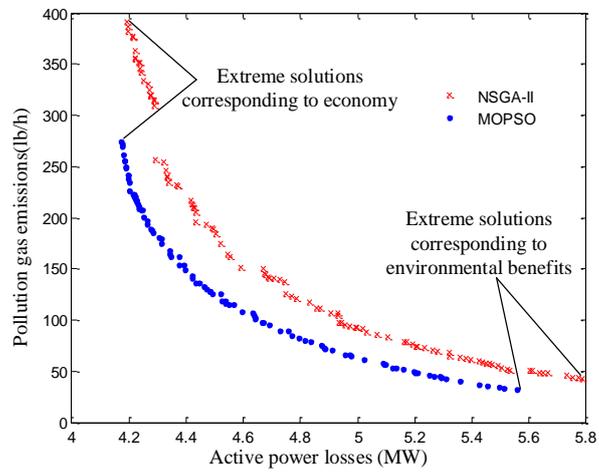

**Fig. 6. Distribution of Pareto frontiers obtained by NSGA-II and MOPSO**

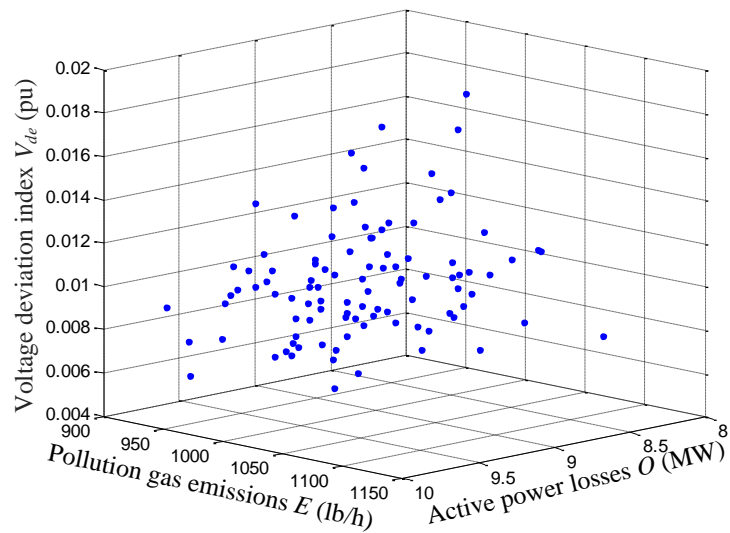

**Fig. 7. Distribution of Pareto-optimal solutions of IEEE 14-bus system with**

**two-terminal DC network**

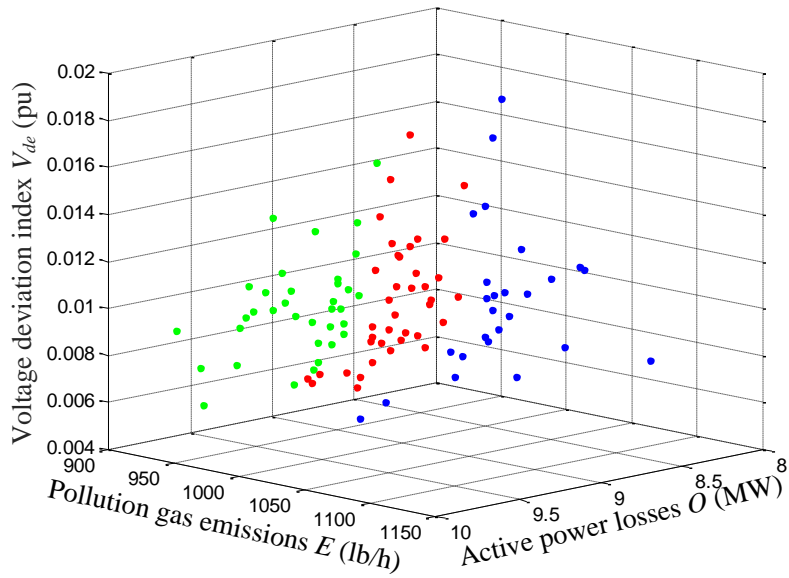

**Fig. 8. Distribution of Pareto-optimal solutions of IEEE 14-bus system with two-terminal DC network after clustering**

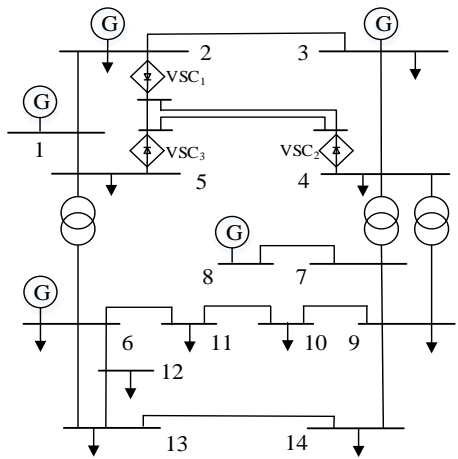

**Fig. 9. IEEE 14-bus system with three-terminal DC network**

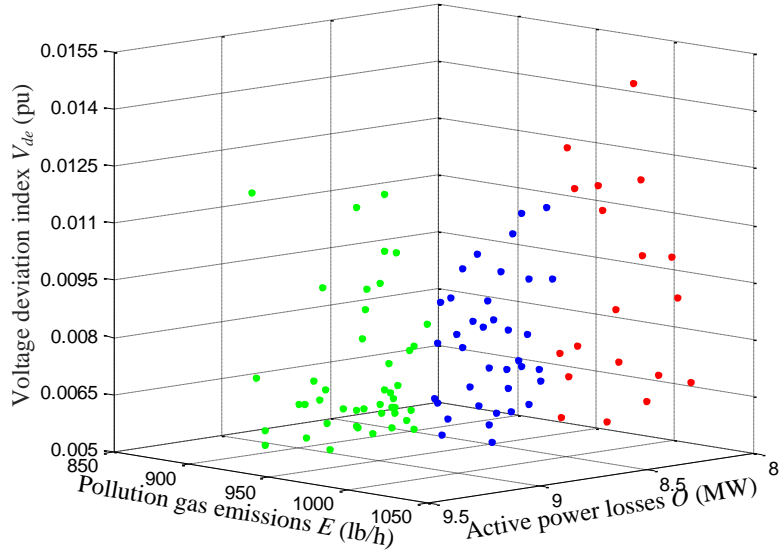

**Fig. 10. Distribution of Pareto-optimal solutions of IEEE 14-bus system with three-terminal DC network**

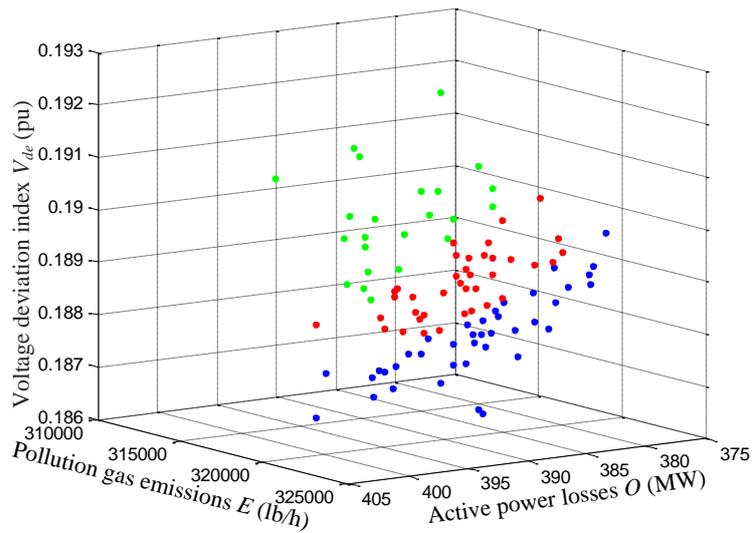

**Fig. 11. Distribution of Pareto-optimal solutions of IEEE 300-bus system**

**Table 1**

Parameter values of MOPSO

| Parameters | Values |
| --- | --- |
| Population size | 100 |
| Repository size | 100 |
| Maximum iteration number | 50 |
| Inertia weight | 0.73 |
| Inertia weight damping ratio | 1 |
| Personal learning coefficient | 1.5 |
| Global learning coefficient | 1.5 |
| Mutation rate | 0.5 |
| Divisions for the adaptive grid | 30 |

**Table 2**

DC bus parameters of the two-terminal DC network

| $B_{VSC}$ | $R$ / pu | $X$ / pu | $P_s$ / pu | $Q_s$ / pu | $U_{dc}$ / pu |
|---|---|---|---|---|---|
| 4 | 0.0015 | 0.1211 | 0.492 | 0.116 | 1.000 |
| 5 | 0.0015 | 0.1211 | -0.495 | -0.105 | 1.000 |

**Table 3**

Extreme solutions in Pareto-optimal solutions

| Algorithm | Optimization objective | $O$ (MW) | $E$ (lb/h) |
|---|---|---|---|
| NSGA-II | Economy | 4.19 | 437 |
| | Environmental benefits | 5.79 | 42 |
| MOPSO | Economy | 4.17 | 275 |
| | Environmental benefits | 5.56 | 32 |

**Table 4**

The maximum, minimum and average number of iterations

| Algorithm | Maximum | Minimum | Average |
|---|---|---|---|
| NSGA-II | 35 | 17 | 23.73 |
| MOPSO | 21 | 12 | 16.53 |

**Table 5**

The maximum, minimum and average times of each algorithm

| Algorithm | Maximum time (s) | Minimum time (s) | Average time(s) |
|---|---|---|---|
| NSGA-II | 27.25 | 35.11 | 31.46 |
| MOPSO | 22.87 | 27.49 | 25.08 |

**Table 6**

Statistical results of GD and SP for the two MOEAs

| Algorithm | Criterions | Best | Average | Worst |
|---|---|---|---|---|
| NSGA-II | GD | 0.1472 | 0.1935 | 0.2404 |
|  | SP | 0.8636 | 1.2723 | 1.6857 |
| MOPSO | GD | 0.1266 | 0.1610 | 0.2072 |
|  | SP | 0.8641 | 1.1330 | 1.3617 |

**Table 7**

Extreme solutions of IEEE 14-bus system with two-terminal DC network

| Extreme solutions | $O$ (MW) | $E$ (lb/h) | $V_{de}$ (pu) |
|---|---|---|---|
| Extreme solution I | 8.18 | 1033.03 | 0.0107 |
| Extreme solution II | 9.59 | 900.54 | 0.0085 |
| Extreme solution III | 9.52 | 912.33 | 0.0054 |

**Table 8**

'Best' compromise solutions of IEEE 14-bus system with 2-terminal DC network

| Compromise solutions | $O$ (MW) | $E$ (lb/h) | $V_{de}$ (pu) | Priority membership |
|---|---|---|---|---|
| Compromise solution I | 8.65 | 935.41 | 0.0102 | 0.7409 |
| Compromise solution II | 9.07 | 903.42 | 0.0095 | 0.6807 |
| Compromise solution III | 8.68 | 998.66 | 0.0063 | 0.7638 |

**Table 9**

VSC variables before and after optimization

| VSC | Before optimization | | | After optimization | | |
|---|---|---|---|---|---|---|
| | $P_s$/pu | $Q_s$/pu | $U_{dc}$/pu | $P_s$/pu | $Q_s$/pu | $U_{dc}$/pu |
| VSC$_1$ | -0.495 | -0.105 | 1.000 | -0.522 | -0.095 | 0.996 |
| VSC$_2$ | 0.492 | 0.116 | 1.000 | 0.495 | 0.104 | 1.049 |

**Table 10**

Generator variables before and after optimization

| Generators | Before optimization | | | After optimization | | |
|---|---|---|---|---|---|---|
| | $P_G$/pu | $Q_G$/pu | $U_G$/pu | $P_G$/pu | $Q_G$/pu | $U_G$/pu |
| G$_1$ | 2.324 | -0.165 | 1.060 | 1.118 | 0.071 | 1.060 |
| G$_2$ | 0.400 | 0.436 | 1.045 | 1.002 | 0.101 | 1.048 |
| G$_3$ | 0 | 0.251 | 1.010 | 0.300 | 0.138 | 1.017 |
| G$_4$ | 0 | 0.127 | 1.070 | 0.099 | 0.197 | 1.035 |
| G$_5$ | 0 | 0.176 | 1.090 | 0.100 | 0.003 | 1.033 |

**Table 11**

Objective functions before and after optimization

| Optimization situation | Network loss rate (%) | $E$ (lb/h) | $V_{de}$ (pu) |
|---|---|---|---|
| Before optimization | 2.09 | 1405.41 | 0.0232 |
| After optimization | 1.12 | 935.41 | 0.0102 |

**Table 12**

DC bus parameters in three-terminal DC network

| $B_{VSC}$ | $R$ / pu | $X$ / pu | $P_s$ / pu | $Q_s$ / pu | $U_{dc}$ / pu |
|---|---|---|---|---|---|
| 2 | 0.006 | 0.150 | -0.839 | 0.142 | 1.000 |
| 4 | 0.006 | 0.150 | 0.968 | 0.016 | 1.000 |
| 5 | 0.006 | 0.150 | -0.129 | 0.134 | 1.000 |

**Table 13**

'Best' compromise solutions of IEEE 14-bus system with 3-terminal DC network

| Compromise solutions | $O$ (MW) | $E$ (lb/h) | $V_{de}$ (pu) | Priority membership |
|---|---|---|---|---|
| Compromise solution I | 8.26 | 973.30 | 0.0075 | 0.7226 |
| Compromise solution II | 8.92 | 903.06 | 0.0092 | 0.6915 |
| Compromise solution III | 8.60 | 931.93 | 0.0061 | 0.6324 |

**Table 14**

VSC variables before and after optimization

| VSCs | Before optimization | | | After optimization | | |
|---|---|---|---|---|---|---|
| | $P_s$/pu | $Q_s$/pu | $U_{dc}$/pu | $P_s$/pu | $Q_s$/pu | $U_{dc}$/pu |
| $VSC_1$ | -0.839 | 0.142 | 1.000 | -0.889 | 0.126 | 1.059 |
| VSC$_2$ | 0.968 | 0.016 | 1.000 | 0.979 | 0.005 | 1.019 |
| VSC$_3$ | -0.129 | 0.134 | 1.000 | -0.143 | 0.117 | 1.040 |

**Table 15**

Generator variables before and after optimization

| Generators | Before optimization | | | After optimization | | |
|---|---|---|---|---|---|---|
| | $P_G$/pu | $Q_G$/pu | $U_G$/pu | $P_G$/pu | $Q_G$/pu | $U_G$/pu |
| G$_1$ | 2.324 | -0.165 | 1.060 | 1.056 | -0.070 | 1.060 |
| G$_2$ | 0.400 | 0.436 | 1.045 | 1.112 | 0.168 | 1.049 |
| G$_3$ | 0 | 0.251 | 1.010 | 0.300 | 0.212 | 1.027 |
| G$_4$ | 0 | 0.127 | 1.070 | 0.099 | 0.083 | 1.029 |
| G$_5$ | 0 | 0.176 | 1.090 | 0.100 | 0.178 | 1.060 |

**Table 16**

Objective functions before and after optimization

| Optimization situation | Network loss rate (%) | $E$ (lb/h) | $V_{de}$ (pu) |
|---|---|---|---|
| Before optimization | 2.09 | 1405.41 | 0.0232 |
| After optimization | 1.07 | 973.30 | 0.0075 |

**Table 17**

Comparative results of different typical operation conditions

| Variables | Maximum | Minimum | Case1 | Case2 | Case3 |
|---|---|---|---|---|---|
| $P_{G1}$/pu | 3.32 | 0.32 | 1.207 | 1.218 | 1.258 |
| $P_{G2}$/pu | 1.40 | 0.40 | 0.972 | 0.964 | 0.924 |
| $P_{G3}$/pu | 0.30 | 0 | 0.300 | 0.299 | 0.300 |
| $P_{G4}$/pu | 0.10 | 0 | 0.099 | 0.100 | 0.100 |
| $P_{G5}$/pu | 0.10 | 0 | 0.100 | 0.097 | 0.098 |
| $U_{G1}$/pu | 1.10 | 0.95 | 1.060 | 1.060 | 1.060 |
| $U_{G2}$/pu | 1.10 | 0.95 | 1.044 | 1.048 | 1.043 |
| $U_{G3}$/pu | 1.10 | 0.95 | 1.011 | 1.012 | 1.007 |
| $U_{G4}$/pu | 1.10 | 0.95 | 1.025 | 1.022 | 1.023 |
| $U_{G5}$/pu | 1.10 | 0.95 | 1.042 | 1.042 | 1.059 |
| $P_{s,1}$/pu | 1.00 | -1.00 | -0.527 | -0.550 | -0.880 |
| $Q_{s,1}$/pu | 1.00 | -1.00 | -0.095 | -0.094 | 0.128 |
| $U_{dc,1}$/pu | 1.10 | 0.95 | 0.988 | 0.988 | 1.052 |
| $P_{s,2}$/pu | 1.00 | -1.00 | 0.501 | 0.523 | 0.969 |
| $Q_{s,2}$/pu | 1.00 | -1.00 | 0.104 | 0.105 | 0.008 |
| $U_{dc,2}$/pu | 1.10 | 0.95 | 1.054 | 1.057 | 1.005 |
| $P_{s,3}$/pu | 1.00 | -1.00 | - | - | -0.143 |
| $Q_{s,3}$/pu | 1.00 | -1.00 | - | - | 0.118 |
| $U_{dc,3}$/pu | 1.10 | 0.95 | - | - | 1.029 |
| $T_1$/pu | 1.1 | 0.9 | 1.0405 | 1.0280 | 1.0155 |
| $T_2$/pu | 1.1 | 0.9 | 0.9565 | 0.9815 | 0.9940 |
| $T_3$/pu | 1.1 | 0.9 | 1.0070 | 1.0070 | 0.9570 |
| $Q_{C1}$/pu | 0.25 | 0 | 0.17 | 0.17 | 0.17 |
| $O$/MW | - | - | 8.87 | 8.92 | 8.80 |
| $E$/(lb/$h$) | - | - | 928.48 | 928.21 | 917.38 |

| | | | 0.0072 | 0.0074 | 0.0053 |
| --- | --- | --- | --- | --- | --- |
| $V_{de}$/pu | - | - | 0.0072 | 0.0074 | 0.0053 |

**Table 18**

'Best' compromise solutions of IEEE 300-bus system

| Compromise solutions | $O$ (MW) | $E$ (lb/h) | $V_{de}$ (pu) | Priority membership |
| --- | --- | --- | --- | --- |
| Compromise solution I | 382.32 | 316016 | 0.1889 | 0.7284 |
| Compromise solution II | 383.19 | 314242 | 0.1897 | 0.5459 |
| Compromise solution III | 386.72 | 319036 | 0.1880 | 0.7194 |

**Table 19**

Comparative results before and after optimization of IEEE 300-bus system

| Optimization situation | Network loss rate (%) | $E$ (lb/h) | $V_{de}$ (pu) |
| --- | --- | --- | --- |
| Before optimization | 1.31 | 340453 | 0.2465 |
| After optimization | 1.17 | 316016 | 0.1889 |